\documentclass[12pt]{article}

\usepackage{amsthm,amsmath,mathrsfs,amsfonts,amssymb,mathtools,nicefrac}
\usepackage{graphicx,caption,subcaption}
\usepackage{enumerate}
\usepackage{natbib} \RequirePackage[colorlinks,citecolor=blue,urlcolor=blue]{hyperref}
\usepackage{url} 

\usepackage[multiple]{footmisc} \usepackage[capposition=top]{floatrow}

\usepackage{adjustbox}
\usepackage{multirow}

\usepackage{amsmath}
\usepackage{accents}
\newlength{\dhatheight}

\usepackage{amsmath}  \usepackage{amsthm}
\usepackage{amsfonts} \usepackage{amssymb}
\usepackage{appendix} \usepackage{graphicx}
\usepackage{mathrsfs} \usepackage{lscape} \usepackage{fixmath} \usepackage{enumitem} \usepackage{tikz} \usepackage{float}
\usepackage[multiple]{footmisc}
\usepackage{,psfrag,epsf,caption,subcaption}
\usepackage[sans]{dsfont}

\usepackage{mathtools}
\newcommand{\egreg}{\vphantom{-}}

\usepackage{nicefrac}

\usetikzlibrary{automata}

\everymath=\expandafter{\the\everymath\textstyle}

\newcommand{\blind}{1}
\newtheorem{assumption}{Assumption}

\newtheorem{proposition}{Proposition}
\newtheorem{theorem}{Theorem}

\theoremstyle{definition}

\begin{document}

\def\spacingset#1{\renewcommand{\baselinestretch}{#1}\small\normalsize} \spacingset{1}

\if1\blind
{
  \title{\bf {\normalsize ESTIMATION AND INFERENCE FOR STOCHASTIC BLOCK MODELS}}
  \author{{\normalsize Koen Jochmans}\thanks{Address: Toulouse School of Economics, 1 esplanade de l'Universit\'e, 31080 Toulouse, France. E-mail: \texttt{koen.jochmans@tse-fr.eu}.} \\  {\small Toulouse School of Economics}, {\small University of Toulouse Capitole}
    }
\date{\small This version: \today}
  \maketitle
} \fi

\if0\blind
{
  \bigskip
  \bigskip
  \bigskip
  \begin{center}
    {\bf {\normalsize ESTIMATION AND INFERENCE FOR STOCHASTIC BLOCK MODELS}}
\end{center}
  \medskip
} \fi

\vspace{-.5cm}
\begin{abstract}
\noindent
This paper is concerned with nonparametric estimation of the weighted stochastic block model. We first show that the model implies a set of multilinear restrictions on the joint distribution of edge weights of certain subgraphs involving (in its simplest form) triplets and quadruples of nodes. From this system of equations the unknown components of the model can be recovered nonparametrically, up to the usual labeling ambiguity. We introduce a simple and computationally-attractive manner to do this. Estimators then follow from the analogy principle. Limit theory is provided. We find that component distributions and their functionals, as well as their density functions (for the case where edge weights are continuous) are all estimable at the parametric rate. Numerical experiments are reported on.

\end{abstract}




\spacingset{1.45}

\section{Introduction}  
The stochastic block model provides a parsimonious way to incorporate latent heterogeneity in the analysis of network data. The original application of  \citealt{HollandLaskeyLeinhardt1983} concerned the binary decision of edge formation between two nodes and generalizes the Erd\H{o}s-R\'enyi random-graph model (\citealt{ErdosRenyi1959}). In the latter model, edges are formed independently with a common probability. In the former, the set of nodes is partitioned into a finite set of latent communities, and the link probability between two nodes depends on the communities that they belong to. The latent block structure has since been used to study general (discrete or continuous) outcomes generated from pairwise interaction (\citealt{HoffRafteryHandcock2002}), thereby extending the applicability of the stochastic block model to weighted graphs.

There is now a large literature on the estimation of the stochastic block model. \cite{Lei2016}, \cite{Wangbickel2017}, \cite{YanSarkarCheng2018}, and \cite{LeLevina2019} provide techniques to estimate and test the number of communities. Taking this number as given, \cite{SnijdersNowicki1997}, \cite{NowickiSnijders2001}, and \cite{AminiChenBickelLevina2013} consider likelihood-based estimation of the remaining parameters of the model. Given the well-known computational complexity of this approach variational methods have also been considered (\citealt{DaudinPicardRobin2008}, \citealt{MariadassouRobinVacher2010}). For the binary edge-weight case statistical properties under growth rates on the average degree have been derived by \cite{CelisseDaudinPierre2012} and \cite{BickelChoiChangZhang2013}. Related results are also available for approaches based on spectral clustering (\citealt{RoheChatterjeeYu2011}, \citealt{SussmanTangFishkindPriebe2012}, \citealt{TangCapePriebe2022}). 

This paper is concerned with estimation and inference for the weighted stochastic block model. We present a set of multilinear equations from which all the unknown parameters---the number of communities, their size, and the conditional distributions of the edge weights---can be uniquely recovered. These restrictions involve the probability distributions of edge weights in small subgraphs of certain configurations. As they hold without imposing any parametric structure, our equations can be used as a basis for the construction of a fully nonparametric estimator of the stochastic block model. It would also be possible to complement our approach with (semi-) parametric restrictions, although we do not focus on this here.

The estimation strategy we lay out below is computationally attractive. It is built around a joint (approximate) diagonalization step tailored to our setup. This type of routine has found applicability elsewhere (\citealt{CardosoSouloumiac1993}, \citealt{BonhommeJochmansRobin2014}). Here we use it as an auxiliary first-step estimator in the construction of our main estimators of the components of the model. Moreover, once it has been computed, our estimators of the stochastic block model are least-squares estimators and are thus immediate to compute. We present an estimator of the distribution of the communities as well as a generic estimator of linear functionals of the conditional distributions. The latter covers (cumulative) distribution functions, their moments, and probability mass functions, for example. We also give results for a kernel estimator for conditional densities for the case where edge weights are continuous. 

Limit theory is presented under an asymptotic scheme where the number of nodes in the network, $n$, goes to infinity, assuming that the number of communities is known. Under weak regularity conditions they converge in distribution to correctly-centered normal random variables at the rate $n^{-\nicefrac{1}{2}}$. Interestingly, this (parametric) rate equally applies to density estimation, as the smoothing bias is small relative to the standard deviation. This is due to the strong dependence induced by the community structure. Undersmoothing is not needed to achieve this result. 

Finally, because our estimators involve averages over only triplets and quadruples of nodes, the conditions underlying our limit results do not impose requirements on network denseness through, for example, (functionals of) the degree distribution. Indeed, they do not attempt to assign nodes to communities. Therefore, our techniques can easily be adapted to a setting where we observe many (possible small) networks generated by the same block structure.

\section{Stochastic block model}
Consider a graph involving $n$ nodes where the set of nodes is partitioned into $r$ latent communities, labelled $1,\ldots, r$. Each node is first assigned to a community independently according to some probability distribution $\boldsymbol{p}=(p_1,\ldots,p_r)^\prime$. The community of node $i$ is recorded in the latent variable $Z_i$. Thus, 
$$
\mathbb{P}(Z_i = z) = p_z>0 \text { for each } 1\leq z \leq r
$$ 
and is equal to zero otherwise. Next, each unordered pair of nodes $i\neq j$ draws a real-valued weighted edge $X_{i,j}$ from some distribution that depends on the communities they belong to, 
$$
F_{z_1,z_2}(x) := \mathbb{P}(X_{i,j}\leq x \vert Z_i=z_1, Z_j=z_2).
$$
The edge weights $X_{i,j}$ are independent conditional on the community indicators $(Z_i,Z_j)$.


\subsection{Multilinear restrictions}
We first show that it is possible to recover the parameters of the stochastic block model from the distribution of edge weights of subgraphs involving as little as four nodes. To do so we let
$$
F_z(x) 
:= 
\mathbb{P}(X_{i,j}\leq x \vert Z_i = z)
=
\sum_{z^\prime=1}^r  F_{z,z^\prime}(x) \, p_{z^\prime }
$$
and impose the following condition.

\begin{assumption} \label{ass:rank}
The functions $F_1,\ldots,F_r$ are linearly independent.
\end{assumption}

\noindent
Rank conditions as this one arise frequently in the analysis of multivariate latent-variable models.

We will prove the following theorem and specialize it further afterwards. 

\begin{figure}[tp] \caption{Subgraph configurations used for identification} \label{fig:subgraphs}
\resizebox{\columnwidth}{!}{	
    \begin{subfigure}[t]{.40\textwidth}
    \caption{$r$} \centering
     \begin{tikzpicture}
     \node[state,fill=white] (v1) at (-1,-1) {$1$};
 	\node[state,fill=white] (v2) at ( 1,-1) {$2$};
 	\node[state,fill=white] (v3) at (-1, 1) {$3$};
 	\draw[blue] (v1) -- (v2);
 	\draw[blue] (v1) -- (v3);
    \end{tikzpicture}
   \end{subfigure} 
   \begin{subfigure}[t]{.40\textwidth}
   \caption{$p_z$ and $F_z$} \centering	
    \begin{tikzpicture}
    \node[state,fill=white] (v1) at (-1,-1) {$1$};
	\node[state,fill=white] (v2) at ( 1,-1) {$2$};
	\node[state,fill=white] (v3) at (-1, 1) {$3$};
	\node[state,fill=white] (v4) at ( 1, 1) {$4$};
	\draw[blue] (v1) -- (v2);
	\draw[blue] (v1) -- (v3);
	\draw[blue] (v1) -- (v4);
   \end{tikzpicture}   
   \end{subfigure}
   \begin{subfigure}[t]{.40\textwidth}
   \caption{$F_{z_1,z_2}$} \centering
    \begin{tikzpicture}
    \node[state,fill=white] (v1) at (-1,-1) {$1$};
	\node[state,fill=white] (v2) at ( 1,-1) {$2$};
	\node[state,fill=white] (v3) at (-1, 1) {$3$};
	\node[state,fill=white] (v4) at ( 1, 1) {$4$};
	\draw[blue] (v1) -- (v2);
	\draw[blue] (v1) -- (v3);
	\draw[blue] (v3) -- (v4);
   \end{tikzpicture}
   \end{subfigure}
   }
\end{figure}
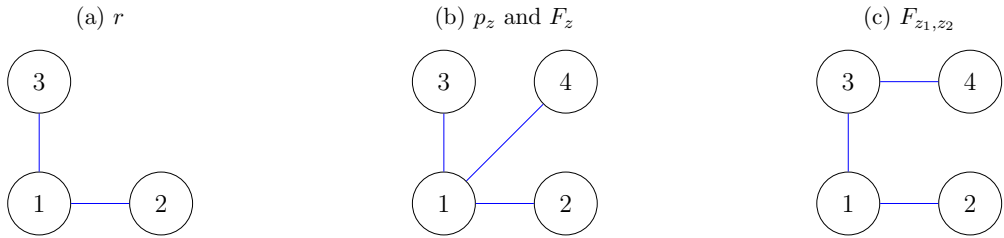

\begin{theorem} \label{thm:id}
Suppose that Assumption \ref{ass:rank} holds. Then,  

\vspace{.2cm}\noindent
(i) \, The number of communities $r$ is nonparametrically recoverable from the distribution of edge weights in two-star subgraphs;

\vspace{.2cm}\noindent
(ii) \, The distribution of communities, $\boldsymbol{p}$, as well as  expectations of the form
$$
\varphi_{z_1,z_2}:=
\mathbb{E}(\varphi(X_{i,j}) \vert Z_i=z_1,Z_j=z_2)
$$
are nonparametrically recoverable (up to the labeling of the latent communities) from the distribution of edge weights in three-star subgraphs and the distribution of edge weights in path subgraphs on four nodes.

\end{theorem}

\noindent
Consider, first, the two-star subgraph in Figure \ref{fig:subgraphs}(a), involving the edges between the three nodes $1,2,3$. Observe that the edge weights $X_{1,2},X_{1,3}$ are independent conditional on $Z_1$. Furthermore, their (unconditional) distribution factors as
$$
\sum_{z=1}^r \, p_z \, F_z \otimes F_z
$$
This is a bivariate finite-mixture model. It follows from the work of \citet[Propositions 2.1 and 2.3]{MbakopKwon2019} that $r$ is identified under Assumption \ref{ass:rank}. See also \cite{KasaharaShimotsu2014} for related results.

Next look at the three-star subgraph in Figure \ref{fig:subgraphs}(b). The edge weights in this three-star graph---$X_{1,2}, X_{1,3}, X_{1,4}$---are again independent conditional on $Z_1$. Hence, their tri-variate distribution again factors as the multivariate finite mixture
$$
\sum_{z=1}^r \, p_z \, F_z \otimes F_z \otimes F_z.
$$
From this, the identification of $F_{1},\ldots, F_{r}$ and  $\boldsymbol{p}=(p_1,\ldots, p_r)^\prime$ (up to an arbitrary but common ordering) follows from \citet[Corollaries 1 and 2]{BonhommeJochmansRobin2014}. See also \cite{AllmanMatiasRhodes2009} for a related result in multivariate mixture models.

Assumption \ref{ass:rank} implies that there exists a finite integer $l$ such that the $l\times r$ matrix $\boldsymbol{G}$,
$$
(\boldsymbol{G})_{l^\prime,z} := \mathbb{E}(\alpha_{l^\prime}(X_{i,j}) \vert Z_i=z),
$$
for a set of transformation functions $\alpha_1,\ldots,\alpha_l$, has full column rank. One choice for these transformation functions would be $\alpha_{l^\prime}(x)= \lbrace x\leq x_{l^\prime} \rbrace$, where  $x_1,\ldots, x_l$ is a grid of points and $\lbrace \cdot \rbrace$ denotes the indicator function. In this case, 
$
(\boldsymbol{G})_{l^\prime,z} = F_z(x_{l^\prime}).
$
Other approximating functions such as orthogonal polynomials are also possible. Whatever the choice of functions, the matrix $\boldsymbol{G}$ is identified because the $F_z$ are and the columns of $\boldsymbol{G}$ are linear functionals thereof. Moving on, use the joint distribution of the edge weights $X_{1,2}, X_{1,3},X_{3,4}$ from the path graph in Figure \ref{fig:subgraphs}(c) to construct the $l \times l$ matrix $\boldsymbol{M}_\varphi$ with
$$
(\boldsymbol{M}_\varphi)_{l_1,l_2}
:=
\mathbb{E}(
\alpha_{l_1}(X_{1,2})
\, 
\varphi(X_{1,3})
\,
\alpha_{l_2}(X_{3,4})
).
$$
Noting that $X_{1,2}, X_{1,3},X_{3,4}$ are independent conditional on the pair $(Z_1,Z_3)$, we have that
$$
\boldsymbol{M}_\varphi
=
\boldsymbol{G}^{\hphantom{\prime}} \boldsymbol{H}_\varphi \boldsymbol{G}^\prime,
$$
where
$$
(\boldsymbol{H}_\varphi)_{z_1,z_2}
:=
p_{z_1} p_{z_2} \varphi_{z_1,z_2}.
$$
Because $\boldsymbol{G}$ has maximal column rank $\boldsymbol{G}^\prime \boldsymbol{G}$ is invertible and, thus, we obtain
$$
\boldsymbol{H}_\varphi
=
(\boldsymbol{G}^\prime \boldsymbol{G})^{-1} \boldsymbol{G}^\prime
\boldsymbol{M}_\varphi \boldsymbol{G} (\boldsymbol{G}^\prime \boldsymbol{G})^{-1}.
$$
With the $p_z$ already shown to be identified this result suffices to show the theorem. It is nonetheless useful to note that an application of the above argument to the constant function $\varphi(x)=1$ gives
$
\boldsymbol{M}_1
=
\boldsymbol{G}^{\hphantom{\prime}} \boldsymbol{H}_1 \boldsymbol{G}^\prime
$
where 
$
(\boldsymbol{H}_1)_{z_1,z_2}
:=
p_{z_1} p_{z_2}.
$
It therefore follows that
$$
\varphi_{z_1,z_2}
=
(\boldsymbol{H}_1)_{z_1,z_2}^{-1} (\boldsymbol{H}_\varphi)_{z_1,z_2}^{\egreg},
$$
which is a convenient result for the purpose of estimation. The proof of Theorem \ref{thm:id} is complete. 

Part (ii) of Theorem \ref{thm:id} can be applied to $\varphi(x)= \lbrace x\leq x^\prime \rbrace$ for any chosen value $x^\prime$, leading to our first proposition.

\begin{proposition} \label{prop:F}
Suppose that Assumption \ref{ass:rank} holds. Then the (weighted) stochastic block model is nonparametrically identified up to relabeling of the latent communities.
\end{proposition}	

\noindent
Proposition \ref{prop:F} is to be contrasted with the existing identification results to date. 
For the unweighted model with two communities \citet[Theorem 7]{AllmanMatiasRhodes2009} showed identification from the (joint distribution of) complete subgraphs involving 16 nodes. Assuming the number of communities is known, \citet[Theorems 14 and 15]{AllmanMatiasRhodes2011} obtained results for the general model, but they rely on the complete subgraph on 9 nodes to do so. Their results further require conditions on the support of the edge weights relative to the number of communities being sufficiently large, in addition to $F_{z_1,z_2}$ (for all $z_1<z_2$) being linearly independent.


The proof of Theorem \ref{thm:id} is remarkably simple. Yet it provides a constructive approach to estimation and we will consider this below. Note, also, that Theorem \ref{thm:id} is useful beyond as an input to establish Proposition \ref{prop:F} as it can be used to directly show identification---and to construct estimators---of linear functionals of the conditional distributions without the need to first estimate the latter.

\subsection{A generalization}
Assumption \ref{ass:rank} cannot be satisfied when the edge weights can take on strictly less values than there are latent communities. 
In such a case identification can be obtained by looking at subgraphs involving a larger set of nodes. This illustrates the interplay between the richness of the support of the edge weights, the number of latent communities, and the size of the subgraphs that are needed to show identification.

The argument is based on chaining nodes in a particular manner. We let
$$
F_z^q(x_1,\ldots,x_q)
:=
\mathbb{P}(X_{i,i_1}\leq x_1, X_{i_1, i_2}\leq x_2,\ldots, X_{i_1,i_{q}} \leq x_q\vert Z_i=z),
$$
where the indices $i$ and $i_1,\ldots,i_q$ are all distinct, and impose the following rank requirement.

\renewcommand{\theassumption}{\arabic{assumption}'}
\setcounter{assumption}{0}

\begin{assumption} \label{ass:rank2}
There exists a finite integer $q$ such that the functions $F_1^q,\ldots,F_r^q$ are linearly independent.
\end{assumption}

\renewcommand{\theassumption}{\arabic{assumption}}
\setcounter{assumption}{1}

\renewcommand{\theproposition}{\arabic{proposition}'}
\setcounter{proposition}{0}

\begin{proposition}\label{prop:F2}
	Let $q_\mathrm{min}$ be the smallest integer for which Assumption \ref{ass:rank2} holds. If $q_\mathrm{min}$ is positive the (weighted) stochastic block model is nonparametrically identified up to relabeling of the latent communities.
\end{proposition}	

\renewcommand{\theproposition}{\arabic{proposition}}
\setcounter{proposition}{1}


\section{Nonparametric estimation}
An estimator of the number of latent communities, $r$, can be constructed along the lines of \cite{MbakopKwon2019}. Here we construct estimators of the distribution of the communities, $\boldsymbol{p}$, and the conditional distributions $F_{z_1,z_2}$ and functionals thereof, building on the proof of Theorem \ref{thm:id}. 
We consider a setting where we observe data from a single network involving $n$ nodes.

Our proposal is to proceed in two sequential steps. First, the matrix $\boldsymbol{G}$ is estimated by a modification of the diagonalization estimator of \cite{BonhommeJochmansRobin2014}. This estimator, $\boldsymbol{\hat{G}}$, is detailed below. Next, we appeal to the analogy principle to construct our estimators of the components of the stochastic block model. An alternative to the joint diagonalization approach would be to estimate the $(\boldsymbol{G})_{l^\prime,z} = \mathbb{E}(\alpha_{l^\prime}(X_{i,j})\vert Z_i=z)$ using estimates of the $F_z$. These could be obtained by maximizing a parametric likelihood (using the EM algorithm, see \citealt{McLachlanPeel2000}), or by nonparametric procedures such as those given in \cite{LevineHunterChauveau2011}. A practical advantage of our proposal is that it bypasses estimation of the complete mixture model. A theoretical advantage (relative to other nonparametric estimators) is that distribution theory for the matrix $\boldsymbol{\hat{G}}$ can be obtained by adapting the work of \cite{BonhommeJochmansRobin2014} to deal with the network structure of the data.

If we write $\boldsymbol{a}:=(a_1,\ldots, a_l)^\prime$ for $a_{l^\prime}:=\mathbb{E}(\alpha_{l^\prime}(X_{i,j}))$ we have the univariate mixture representation 
$$
\boldsymbol{a}=\boldsymbol{G}\boldsymbol{p}.
$$
Given $\boldsymbol{\hat{G}}$, a least-squares argument suggests estimating $\boldsymbol{p}$ by
$$
\boldsymbol{\hat{p}} := (\boldsymbol{\hat{G}}^\prime \boldsymbol{\hat{G}})^{-1} \boldsymbol{\hat{G}}^\prime \boldsymbol{\hat{a}},
$$ 
for $\boldsymbol{\hat{a}}=(\hat{a}_1,\ldots, \hat{a}_l)^\prime$, with
$$
\hat{a}_{l^\prime}:= \frac{2}{n(n-1)}
\sum_{i < j} \alpha_{l^\prime}(X_{i,j}).
$$
This approach is inspired by \cite{Titterington1983}, where minimum-distance estimators of mixing proportions were considered.

Similarly, 
$$
\boldsymbol{\hat{H}}_\varphi
:=
(\boldsymbol{\hat{G}}^\prime \boldsymbol{\hat{G}})^{-1} \boldsymbol{\hat{G}}^\prime \boldsymbol{\hat{M}}_\varphi
\boldsymbol{\hat{G}}
(\boldsymbol{\hat{G}}^\prime \boldsymbol{\hat{G}})^{-1},
$$
constructed with
$$
(\boldsymbol{\hat{M}}_\varphi)_{l_1,l_2}
:=
\frac{1}{n(n-1)(n-2)(n-3)}
\hspace{-.1cm}
\sum_{i_1\neq i_2\neq i_3 \neq i_4} \hspace{-.3cm}
\alpha_{l_1}(X_{i_1,i_2})
\, 
\varphi(X_{i_2,i_3})
\,
\alpha_{l_2}(X_{i_3,i_4}),
$$
yields the estimator
$$
\hat{\varphi}_{z_1,z_2} :=
(\boldsymbol{\hat{H}}_1)_{z_1,z_2}^{-1} \, (\boldsymbol{\hat{H}}_\varphi)_{z_1,z_2}^{\egreg}
$$
of $\varphi_{z_1,z_2}$. 

Below we will present the sampling properties of these estimators under asymptotics where the number of nodes, $n$, grows large, assuming that the number of communities, $r$, is known. 
An alternative sampling scheme would be to sample $m$ independent networks, each of size $n$ and generated from the same stochastic block model. Under asymptotics where $m \rightarrow \infty$ while $n$ remains fixed, our estimators achieve the parametric rate of $m^{-\nicefrac{1}{2}}$ under the same regularity conditions as the ones introduced here.
We omit further details for this case for brevity.

\subsection{Diagonalization step} 
We construct the matrix $\boldsymbol{\hat{G}}$ by relying on an (approximate) simultaneous-diagonalization argument. We summarize the procedure here and refer to \cite{BonhommeJochmansRobin2014} for additional details on this approach in (stationary) multivariate mixture models. We begin by constructing the $l\times l$ matrix $\boldsymbol{\hat{A}}_0$, 
$$
(\boldsymbol{\hat{A}}_0)_{l_1,l_2}
:=
\frac{1}{n(n-1)(n-2)}
\sum_{i_1\neq i_2\neq i_3} \alpha_{l_1}(X_{i_1,i_2}) \, \alpha_{l_2}(X_{i_1,i_3}),
$$
and perform an eigendecomposition on it to construct an $l\times r$ matrix $\boldsymbol{\hat{V}}$ for which $\boldsymbol{\hat{V}} \boldsymbol{\hat{A}}_0  \boldsymbol{\hat{V}}^\prime = \boldsymbol{I}_r$, the $r\times r$ identity matrix. We next form the $l\times l$ matrices $\boldsymbol{\hat{A}}_{l^\prime}$, 
$$
(\boldsymbol{\hat{A}}_{l^\prime})_{l_1,l_2}
:=
\frac{1}{n(n-1)(n-2)(n-3)}
\hspace{-.1cm}
\sum_{i_1\neq i_2\neq i_3 \neq i_4} \hspace{-.3cm}
\alpha_{l_1}(X_{i_1,i_2}) \, \alpha_{l^\prime}(X_{i_1,i_3}) \,  \alpha_{l_2}(X_{i_1,i_4}),
$$
where $l^\prime=1,\ldots,l$, and transform them using $\boldsymbol{\hat{V}}$ to obtain the $r\times r$ matrices
$$
\boldsymbol{\hat{N}}_{l^\prime}:= 
\boldsymbol{\hat{V}}^\prime
\boldsymbol{\hat{A}}_{l^\prime}
\boldsymbol{\hat{V}}^\prime.
$$
We then find the matrix of joint (approximate) eigenvectors of these matrices as
$$
\boldsymbol{\hat{Q}}:=
\arg\min_{\boldsymbol{Q}\in \boldsymbol{\mathcal{Q}}}
\sum_{l^\prime=1}^l 
\sum_{z_1\neq z_2} (\boldsymbol{Q}^\prime \boldsymbol{\hat{N}}_{l^\prime}\boldsymbol{Q})_{z_1,z_2}^2,
$$
where we let $\boldsymbol{\mathcal{Q}}$ be the set of $r\times r$ orthonormal matrices. With this matrix at hand we construct $\boldsymbol{\hat{G}}$ as
$$
(\boldsymbol{\hat{G}})_{l^\prime,z} :=
(\boldsymbol{\hat{Q}}^\prime \boldsymbol{\hat{N}}_{l^\prime}\boldsymbol{\hat{Q}})_{z,z}.
$$
The minimization problem that defines $\boldsymbol{\hat{Q}}$ can be solved efficiently using the algorithm of \cite{CardosoSouloumiac1993}.

Before proceeding to the asymptotic properties of this procedure it is useful to comment on its computational complexity. The procedure relies on estimated matrices that take the form of U-statistics up to order four. We note, however, that the kernels of these U-statistics are multiplicatively separable in the indices. Consequently, with some re-arrangement, their computational complexity is of the same order as that of a sample mean. The supplement provides additional details on this.

\subsection{Regularity conditions} 
Three regularity conditions will be used. They are collected here. The first two of them impose conventional requirements on second moments.

\begin{assumption} \label{ass:moments}
The variables $\alpha_{l^\prime}(X_{i,j})$ have finite variance.
\end{assumption}

\noindent
Note that this assumption can always be satisfied by working with bounded functions.

\begin{assumption} \label{ass:moments2}
The variable $\varphi(X_{i,j})$ has finite variance.
\end{assumption}


The third regularity condition concerns the $l\times l$ matrix $\boldsymbol{A}_0$, with
$$
(\boldsymbol{A}_0)_{l_1,l_2}:= \mathbb{E}(\alpha_{l_1}(X_{i_1,i_2}) \, \alpha_{l_2}(X_{i_1,i_3})).
$$
Observe that $\boldsymbol{A}_0 = \boldsymbol{G}\, \mathrm{diag}(\boldsymbol{p})\, \boldsymbol{G}^\prime$. Its rank is equal to $r$ by Assumption \ref{ass:rank} and so it has $r$ non-zero eigenvalues. We represent its eigendecomposition as   
$$
\boldsymbol{A}_0 = \boldsymbol{U}\boldsymbol{L}\, \boldsymbol{U}^\prime,
$$ 
with $\boldsymbol{L}$ the $r\times r$ diagonal matrix that collects the non-zero eigenvalues and $\boldsymbol{U}$ the $l\times r$ orthonormal matrix whose $r$ columns contain the associated eigenvectors. Then
$$
\boldsymbol{V} := \boldsymbol{L}^{-\nicefrac{1}{2}}\boldsymbol{U}^{\prime}
$$
is the probability limit of $\boldsymbol{\hat{V}}$. 


\begin{assumption} \label{ass:A0}
All non-zero eigenvalues of  $\boldsymbol{A}_0$ are simple.
\end{assumption}	

\noindent
This assumption implies continuity of $\boldsymbol{V}$ as a function of $\boldsymbol{A}_0$ and is helpful in deriving the properties of $\boldsymbol{\hat{V}}$.

We remark that, because $\boldsymbol{\hat{A}}_0$ is a $\sqrt{n}$-consistent estimator of $\boldsymbol{A}_0$, the rank condition on the matrix $\boldsymbol{G}$ can be tested by any of a number of standard procedures to test the rank of a matrix.


\subsection{A linearization}
An important step in deriving the large-sample properties of our procedures lies in analyzing the first-step estimator, $\boldsymbol{\hat{G}}$. This estimator is a complicated function of the auxiliary estimators of the matrices $\boldsymbol{A}_0$ and $\boldsymbol{A}:=(\boldsymbol{A}_1,\ldots,\boldsymbol{A}_l)$, with the elements of the latter equal to
$$
(\boldsymbol{A}_{l^\prime})_{l_1,l_2}:=
\mathbb{E}(\alpha_{l_1}(X_{i_1,i_2}) \, \alpha_{l^\prime}(X_{i_1,i_3}) \, \alpha_{l_2}(X_{i_1,i_4})).
$$
Their respective influence functions are
$$
\boldsymbol{\beta}_i(\boldsymbol{A}_0):=\mathrm{vec}\, (\boldsymbol{B}_{i}(\boldsymbol{A}_0) - \mathbb{E}(\boldsymbol{B}_{i}(\boldsymbol{A}_0))),
$$
where
\begin{equation*}
\begin{split}
(\boldsymbol{B}_{i_1}(\boldsymbol{A}_0))_{l_1,l_2}
 := &
\mathbb{E}(\alpha_{l_1}(X_{i_1,i_2}) \, \alpha_{l_2}(X_{i_1,i_3}) \vert Z_{i_1})
\\
 + &
\mathbb{E}(\alpha_{l_1}(X_{i_1,i_2}) \, \alpha_{l_2}(X_{i_2,i_3}) \vert Z_{i_1})
\\
 + &
\mathbb{E}(\alpha_{l_1}(X_{i_2,i_3}) \, \alpha_{l_2}(X_{i_1,i_2}) \vert Z_{i_1}),
\end{split}
\end{equation*}		
and 
$$
\boldsymbol{\beta}_i(\boldsymbol{A}):=\mathrm{vec}\, (\boldsymbol{B}_{i}(\boldsymbol{A})-\mathbb{E}(\boldsymbol{B}_{i}(\boldsymbol{A}))),
$$
where 
$
\boldsymbol{B}_{i}(\boldsymbol{A}):=
(\boldsymbol{B}_{i}(\boldsymbol{A}_{1}),\ldots, \boldsymbol{B}_{i}(\boldsymbol{A}_{l}))
$
for
\begin{equation*}
\begin{split}
(\boldsymbol{B}_{i_1}(\boldsymbol{A}_{l^\prime}))_{l_1,l_2}
:= &
\mathbb{E}(\alpha_{l_1}(X_{i_1,i_2}) \, \alpha_{l_2}(X_{i_1,i_3}) \, \alpha_{l^\prime}(X_{i_1,i_4}) \vert Z_{i_1})
\\
 + &
\mathbb{E}(\alpha_{l_1}(X_{i_1,i_2}) \, \alpha_{l_2}(X_{i_2,i_3}) \, \alpha_{l^\prime}(X_{i_2,i_4}) \vert Z_{i_1})
\\
 + &
\mathbb{E}(\alpha_{l_1}(X_{i_2,i_3}) \, \alpha_{l_2}(X_{i_1,i_2}) \, \alpha_{l^\prime}(X_{i_2,i_4}) \vert Z_{i_1})
\\
 + &
\mathbb{E}(\alpha_{l_1}(X_{i_2,i_3}) \, \alpha_{l_2}(X_{i_2,i_4}) \, \alpha_{l^\prime}(X_{i_1,i_2}) \vert Z_{i_1}).
\end{split}
\end{equation*}	
From this a linearization of $\boldsymbol{\hat{G}}$ can be derived. Additional notation is needed in order to state the result.

The estimator $\boldsymbol{\hat{G}}$ is due to \cite{CardosoSouloumiac1993} and is based on the insight that the $r\times r$ matrices 
$
\boldsymbol{N}_{l^\prime}:=\boldsymbol{V}\boldsymbol{A}_{l^\prime}\boldsymbol{V}^\prime
$
are diagonalizable in the same (orthonormal) basis. We write $\boldsymbol{Q}$ for the $r\times r$ matrix of joint eigenvectors and let 
$$
\boldsymbol{D}_{l^\prime} := \boldsymbol{Q}^\prime
\boldsymbol{N}_{l^\prime}
\boldsymbol{Q}
$$
be the $r\times r$ diagonal matrices that contain their respective eigenvalues. We note that the main diagonal of $\boldsymbol{D}_{l^\prime}$ corresponds to the $l^\prime$-th row of matrix $\boldsymbol{G}$. If
we let
$\boldsymbol{D}:=(\boldsymbol{D}_1,\ldots,\boldsymbol{D}_l)^\prime$, then
$$
\mathrm{vec} (\boldsymbol{G}^\prime)
=
(\boldsymbol{I}_l \otimes \boldsymbol{S}) \, \mathrm{vec} (\boldsymbol{D}^\prime),
$$
where $\boldsymbol{S}:=(\boldsymbol{s}_1,\ldots,\boldsymbol{s}_r)^\prime$ is the $r\times r^2$ selection matrix with $\boldsymbol{s}_z$ the selection vector whose $((z-1)(r+1)+1)$-th entry is equal to one and all other entries are equal to zero. The shorthand
$
\boldsymbol{R}:=(\boldsymbol{D}_1\ominus \boldsymbol{D}_1,\ldots, \boldsymbol{D}_l\ominus \boldsymbol{D}_l)^\prime,
$
where $\boldsymbol{D}_{l^\prime}\ominus \boldsymbol{D}_{l^\prime}:=(\boldsymbol{D}_{l^\prime} \otimes \boldsymbol{I}_r)-(\boldsymbol{I}_r \otimes \boldsymbol{D}_{l^\prime})$ is a Kronecker difference, will also be useful.

We have
$$
\mathrm{vec}(\boldsymbol{\hat{G}}^\prime - \boldsymbol{\vphantom{\hat{G}}G}^\prime)
=
\frac{1}{n} \sum_{i=1}^n 
\boldsymbol{\beta}_i(\boldsymbol{G}^\prime)
+ o_p(n^{-\nicefrac{1}{2}}),
$$
for
$$
\boldsymbol{\beta}_i(\boldsymbol{G}^\prime)
:=
(\boldsymbol{I}_l \otimes \boldsymbol{S}) \,
((\boldsymbol{I}_l \otimes \boldsymbol{I}_{r^2}) + \boldsymbol{P}_{\boldsymbol{R}})
( (\boldsymbol{I}_l \otimes \boldsymbol{K}_{r^2}) 
(\boldsymbol{T}_1+\boldsymbol{T}_2) \boldsymbol{\beta}_i(\boldsymbol{A}_0)
+
\boldsymbol{T}_3 \, \boldsymbol{\beta}_i(\boldsymbol{A})
),
$$
where $\boldsymbol{P}_{\boldsymbol{R}}:=\boldsymbol{R}(\boldsymbol{R}^\prime \boldsymbol{R})^* \boldsymbol{R}^\prime$ with a $*$ superscript denoting the Moore-Penrose pseudo inverse of a matrix, $\boldsymbol{K}_{r^2}:=\boldsymbol{I}_{r^2}+\boldsymbol{C}_{r^2}$ with $\boldsymbol{C}_{r^2}$ denoting the $r^2\times r^2$ commutation matrix, and
$$
\boldsymbol{T}_1:=
-
(\boldsymbol{D}\otimes \boldsymbol{I}_r)
(\boldsymbol{Q}^\prime \otimes \boldsymbol{Q}^\prime)
(\boldsymbol{L}\ominus \boldsymbol{L})^*
(\boldsymbol{L}\otimes \boldsymbol{I}_r)
(\boldsymbol{V}\otimes \boldsymbol{V}),
$$
and
$$
\boldsymbol{T}_2:=
-\frac{1}{2}
(\boldsymbol{D}\otimes \boldsymbol{I}_r)
(\boldsymbol{Q}^\prime \overset{c}{\otimes} \boldsymbol{Q}^\prime)
(\boldsymbol{V}\overset{r}{\otimes} \boldsymbol{V}),
\qquad
\boldsymbol{T}_3:=
\boldsymbol{I}_l \otimes
(\boldsymbol{Q}^\prime \otimes \boldsymbol{Q}^\prime)
(\boldsymbol{V} \otimes \boldsymbol{V}).
$$
Here, we use $\overset{c}{\otimes}$ and $\overset{r}{\otimes}$ to denote columnwise and rowwise Kronecker products, respectively. 

\subsection{Limit behavior}
Given the large-sample behavior of the first-step estimator it is readily established that
$$
\mathrm{vec} (\boldsymbol{\hat{G}}^* - \boldsymbol{\vphantom{\hat{G}}G}^*)
 =
 \frac{1}{n} \sum_{i=1}^n
\boldsymbol{\beta}_i(\boldsymbol{G}^*) + o_p(n^{-\nicefrac{1}{2}}),
$$	
for
$$
\boldsymbol{\beta}_i(\boldsymbol{G}^*)
:=
 (
((\boldsymbol{I}_l - \boldsymbol{G}\boldsymbol{G}^*) \otimes (\boldsymbol{G}^\prime \boldsymbol{G})^{-1}) 
-
(\boldsymbol{G}^{\prime} \otimes \boldsymbol{G})^{*} \, \boldsymbol{C}_{rl}
)
\,
\boldsymbol{\beta}_i(\boldsymbol{G}^\prime).
$$
We can now present the asymptotic behavior of our estimators of the main  components of the stochastic block model. We state these in the form of two theorems.

We first provide the limit distribution of the estimator $\boldsymbol{\hat{p}}$. The asymptotic variance of this estimator is equal to
$$
\boldsymbol{V}_{\boldsymbol{p}}: = \mathbb{E}(\boldsymbol{\theta}_i\boldsymbol{\theta}_i^\prime),
$$
where
$$
\boldsymbol{\theta}_i:=
\boldsymbol{G}^* \, \boldsymbol{\beta}_i(\boldsymbol{a}) + (\boldsymbol{a}^\prime \otimes \boldsymbol{I}_r) \, \boldsymbol{\beta}_i(\boldsymbol{G}^*),
$$
with $\boldsymbol{\beta}_i(\boldsymbol{a}) := (\beta_i(a_1),\ldots,\beta_i(a_l))^\prime$ for
$$
\beta_i(a_{l^\prime}):= 2 \, (\mathbb{E}(\alpha_{l^\prime}(X_{i,j}) \vert Z_i) - a_{l^\prime}).
$$
Theorem \ref{thm:p} follows.

\begin{theorem} \label{thm:p}
Suppose that Assumptions \ref{ass:rank}, \ref{ass:moments}, and \ref{ass:A0} hold. Then
$$
\sqrt{n}\, 
(\boldsymbol{\hat{p}} - \boldsymbol{p})
 \rightsquigarrow N(\boldsymbol{0}, \boldsymbol{V}_{\boldsymbol{p}}),
$$	
as $n\rightarrow\infty$.
\end{theorem}

We next state the limit distribution of our estimator of $\varphi_{z_1,z_2}$. 
We will need
$$
\boldsymbol{\beta}_i(\boldsymbol{M}_\varphi):=\mathrm{vec}\, (\boldsymbol{B}_{i}(\boldsymbol{M}_\varphi)-\mathbb{E}(\boldsymbol{B}_{i}(\boldsymbol{M}_\varphi))),
$$
with
\begin{equation*}
\begin{split}
(\boldsymbol{B}_{i_1}(\boldsymbol{M}_{\varphi}))_{l_1,l_2}
:= &
\mathbb{E}(\alpha_{l_1}(X_{i_1,i_3}) \, \varphi(X_{i_1,i_2}) \, \alpha_{l_2}(X_{i_2,i_4})  \vert Z_{i_1})
\\
 + &
\mathbb{E}(\alpha_{l_1}(X_{i_2,i_3}) \, \varphi(X_{i_1,i_2}) \, \alpha_{l_2}(X_{i_1,i_4})  \vert Z_{i_1})
\\
 + &
\mathbb{E}(\alpha_{l_1}(X_{i_1,i_2}) \, \varphi(X_{i_2,i_3}) \, \alpha_{l_2}(X_{i_3,i_4})  \vert Z_{i_1})
\\
 + &
\mathbb{E}(\alpha_{l_1}(X_{i_3,i_4}) \, \varphi(X_{i_2,i_3}) \, \alpha_{l_2}(X_{i_1,i_2})  \vert Z_{i_1}) .
\end{split}
\end{equation*}	
We then have that
$$
\boldsymbol{\beta}_i(\boldsymbol{H}_\varphi)
:=
\boldsymbol{K}_{r^2} \,
(\boldsymbol{G}^*\boldsymbol{M}_\varphi \otimes \boldsymbol{I}_r) \,
\boldsymbol{\beta}_i(\boldsymbol{G}^*)
+
(\boldsymbol{G}^* \otimes \boldsymbol{G}^*) \, 
\boldsymbol{\beta}_i(\boldsymbol{M}_\varphi)
$$
is the influence function of $\boldsymbol{\hat{H}}_\varphi$. From this, the asymptotic behavior of $\hat{\varphi}_{z_1,z_2}$ will follow after a linearization. Writing $\boldsymbol{I}_r = (\boldsymbol{e}_1,\ldots, \boldsymbol{e}_r)$, the influence function of $
\hat{\varphi}_{z_1,z_2}$ is
$$
\vartheta_i({z_1,z_2}) :=
(p_{z_1} \, p_{z_2})^{-1}
(\boldsymbol{e}_{z_2}^\prime \otimes \, \boldsymbol{e}_{z_1}^\prime) \, 
(\boldsymbol{\beta}_i(\boldsymbol{H}_\varphi)-
\varphi_{z_1,z_2} \, 
\boldsymbol{\beta}_i(\boldsymbol{H}_1)).
$$
We let
$$
v_\varphi(z_1,z_2) : = \mathbb{E}(\vartheta_i({z_1,z_2})^2)
$$
in the next theorem.

\begin{theorem} \label{thm:varphi}
Suppose that Assumptions \ref{ass:rank}, \ref{ass:moments}, \ref{ass:moments2}, and \ref{ass:A0} hold. Then, 
$$
\sqrt{n} \,
(\hat{\varphi}_{z_1,z_2} - {\varphi}_{z_1,z_2})
\rightsquigarrow N(0,v_\varphi(z_1,z_2)),
$$	
as $n\rightarrow\infty$.
\end{theorem}	

\noindent
This theorem covers distribution functions, probability mass functions when edge weights are discrete, and moments, for example. 

More generally, the result can be used in combination with standard asymptotic theory to construct estimators of a parameter 
defined as
$$
\arg\max_{\boldsymbol{\delta}}\mathbb{E}(\varphi(X_{i,j};\boldsymbol{\delta}) \vert Z_i=z_1, Z_j=z_2).
$$
Under regularity conditions the implied estimator will be $\sqrt{n}$-consistent and asymptotically normal. 

\subsection{Density estimation}
Next consider the case where the edge weights are continuous and $F_{z_1,z_2}$ admits a density function, $f_{z_1,z_2}$, say. The estimation of $f_{z_1,z_2}$ may be of interest. Theorem \ref{thm:varphi} does not immediately cover this as such nonparametric estimators involve smoothing- or truncation bias. However, the dependence between edge weights reduces these issues to second-order problems. This is in line with the conclusion reached in recent related work by \cite{GrahamNiuPowell2019}, who considered estimation of the marginal density of the edge weights in a setting that encompasses ours.

We consider a standard kernel estimator at a point $x$. Similar results to the one to follow may be established for locally-linear (or polynomial) versions of the kernel estimator, as well as for estimators based on series expansions. The kernel estimator can be cast into our generic formula for $\hat{\varphi}_{z_1,z_2}$ by setting
$$
\varphi(X_{i,j}) = 
\frac{1}{h_n} k\left(\frac{X_{i,j}-x}{h_n}\right),
$$
where $k$ is a kernel function and $h_n$ is a non-negative bandwidth.

The following conditions are standard in density estimation.

\renewcommand{\theassumption}{\arabic{assumption}'}
\setcounter{assumption}{2}

\begin{assumption} \label{ass:kernel}
The kernel function $k$ is symmetric, bounded, and integrates to one. 
The density functions $f_{z_1,z_2}$ are bounded and are twice differentiable with bounded derivatives.	
\end{assumption}	

\renewcommand{\theassumption}{\arabic{assumption}}
\setcounter{assumption}{4}

\noindent
In the sequel, Assumption \ref{ass:kernel} replaces Assumption \ref{ass:moments2}.

The dependence that the stochastic block model induces between the edge weights means that the variance of the kernel estimator will be of the order
$$
n^{-1} + (n^2h_n)^{-1} + (n^3h_n)^{-1} + (n^4h_n)^{-1}.
$$
These terms arise from the covariances between the (symmetrized) kernel of
$$
\frac{1}{n(n-1)(n-2)(n-3)}
\hspace{-.1cm}
\sum_{i_1\neq i_2\neq i_3 \neq i_4} \hspace{-.3cm}
\alpha_{l_1}(X_{i_1,i_2})
\, 
\frac{1}{h_n} k\left(\frac{X_{i_2,i_3}-x}{h_n}\right)
\,
\alpha_{l_2}(X_{i_3,i_4}),
$$
evaluated at two different quadruples of nodes
that have exactly one, two, three, or all four indices in common. 
Hence, if $n h_n \rightarrow \infty$, the variance is of order $n^{-1}$, and only terms involving quadruples of nodes that have one index in common contribute to the asymptotic variance. 
By the usual arguments for U-statistics on graphs \citep{JansonNowicki1991} their contribution is equal to the variance of the sample mean of the projection of the kernel onto the $Z_i$. 

Further, exploiting the fact that the three terms in the kernel above are independent conditional on $(Z_{i_2},Z_{i_3})$ it is readily confirmed that a standard argument, as validated by Assumption \ref{ass:kernel}, implies that the expectation of the statistic given in the previous display is 
$$
\sum_{z_1,z_2}(\boldsymbol{G})_{l_1,z_1} \, f_{z_1,z_2}(x) \, p_{z_1} p_{z_2} \, (\boldsymbol{G})_{l_2,z_2} + O(h_n^2).
$$
Thus, the smoothing bias of the kernel estimator is of the order $h^2_n$, which is the conventional result for such a procedure. In light of the variance being of the order $n^{-1}$, this means that asymptotic bias will be absent provided that $nh_n^4\rightarrow 0$. Taking these observations together leads to the following conclusion.

\begin{proposition} \label{prop:density}
Suppose that Assumptions \ref{ass:rank}, \ref{ass:moments}, \ref{ass:kernel}, and \ref{ass:A0} hold. Then, if the bandwidth satisfies $n h_n \rightarrow \infty$ and $nh_n^4\rightarrow 0$ as $n\rightarrow\infty$, Theorem \ref{thm:varphi} applies to the kernel density estimator.
\end{proposition}

\noindent
Note that, from above, the mean-squared error of the density estimator is of the order
$$
h_n^4 + n^{-1} + (n^2 h_n)^{-1}.
$$
Equating the rate of the first term to the rate of the second term gives us the optimal-rate requirement $h_n \propto n^{-2/5}$. Observe that this requirement is compatible with the condition $nh_n^4 \rightarrow 0$ in Proposition \ref{prop:density}. Hence, undersmoothing is not needed to prevent asymptotic bias.

\section{Numerical experiments}
We now provide simulation evidence for the standard block model with binary outcomes and two latent communities. Here, the conditional distributions are fully characterized by their success probabilities and so we consider estimation of
$$
\varphi_{z_1,z_2} =
\mathbb{P}(X_{i,j}=1 \vert Z_i=z_1,Z_j=z_2),
$$
along with the relative sizes of the two latent communities, $p_1$ and $p_2$. We report results for several combinations of these probabilities. For each, we simulated $10,000$ networks of size $n=100$, with $p_2=.70$, and report the mean, median, standard deviation, and interquartile range across the Monte Carlo replications. To give a sense of the numerical complexity, estimation of the model for a single replication takes just under $\nicefrac{1}{3}$ of a second on my desktop computer.

Table \ref{table:mc1} contains results for three designs that feature complementarity, i.e.,  success is more likely if agents are from the same community. The three designs vary in how much $\varphi_{1,1}$ is separated from $\varphi_{2,2}$. The specification is peculiar in that agents from different communities never generate successes. We do this to highlight that such a degeneracy does not cause problems for our procedure. 

The table shows good performance of our procedure. The conditional distributions are accurately recovered. As $\varphi_{2,2}$ moves further away from $\varphi_{1,1}$ the standard deviation of the estimated success probabilities goes down, as expected. The estimator of the  population shares of the communities equally does well across the designs. Its performance is essentially unaffected by the design changes.

\begin{table} 
	\caption{Simulation results}
\begin{tabular}{lccccc} 
	\hline\hline
&	$\varphi_{1,1}$ &	$\varphi_{1,2}$ &$\varphi_{2,2}$ &	$p_1$ &	$p_2$ \\
	\hline
	\multicolumn{6}{c}{Design 1} \\
true value  &	0.200 &	0.000 &	0.400 &	0.300 &	0.700  \\

mean   &	0.220 &	0.000 &	0.392 &	0.285 &	0.715  \\
median &	0.219 &	0.000 &	0.392 &	0.284 &	0.716  \\
std.~dev.    &	0.029 &	0.006 &	0.014 &	0.044 &	0.044  \\
iqr    &	0.039 &	0.008 &	0.018 &	0.059 &	0.059  \\
					\hline
						\multicolumn{6}{c}{Design 2} \\
true value &	0.200 &	0.000 &	0.600 &	0.300 &	0.700 \\

mean   &	0.209 &	0.000 &	0.590 &	0.287 &	0.713 \\
median &	0.209 &	0.000 &	0.591 &	0.287 &	0.714 \\
std.~dev.    &	0.024 &	0.003 &	0.012 &	0.043 &	0.043 \\
iqr    &	0.032 &	0.004 &	0.016 &	0.059 &	0.059 \\
					\hline
					\multicolumn{6}{c}{Design 3} \\	
true value  &	0.200 &	0.000 &	0.800 &	0.300 &	0.700 \\

mean   &	0.202 &	0.000 &	0.789 &	0.287 &	0.713 \\
median &	0.202 &	0.000 &	0.789 &	0.287 &	0.713 \\
std.~dev.    &	0.023 &	0.002 &	0.009 &	0.044 &	0.044 \\
iqr    &	0.030 &	0.002 &	0.012 &	0.058 &	0.058	\\
\hline\hline
\end{tabular}
\label{table:mc1}
\end{table}



\section*{Supplementary material}
The proofs of all the technical results are available in the supplement to this paper \citep{Jochmans2021}.


\small
\setlength{\bibsep}{0pt} 
\bibliographystyle{chicago3} 
\bibliography{bibliography}

\begin{thebibliography}{}

\bibitem[\protect\citeauthoryear{Allman, Matias, and Rhodes}{Allman, Matias and
  Rhodes}{2009}]{AllmanMatiasRhodes2009}
Allman, E.~S., C.~Matias, and J.~A. Rhodes (2009).
\newblock Identifiability of parameters in latent structure models with many
  observed variables.
\newblock {\em Annals of Statistics\/}~{\em 37}, 3099--3132.

\bibitem[\protect\citeauthoryear{Allman, Matias, and Rhodes}{Allman, Matias and
  Rhodes}{2011}]{AllmanMatiasRhodes2011}
Allman, E.~S., C.~Matias, and J.~A. Rhodes (2011).
\newblock Parameter identif\mbox{}iability in a class of random graph mixture
  models.
\newblock {\em Journal of Statistical Planning and Inference\/}~{\em 141},
  1719--1736.

\bibitem[\protect\citeauthoryear{Amini, Chen, Bickel, and Levina}{Amini, Chen,
  Bickel and Levina}{2013}]{AminiChenBickelLevina2013}
Amini, A.~A., A.~Chen, P.~J. Bickel, and E.~Levina (2013).
\newblock Pseudo-likelihood methods for community detection in large sparse
  networks.
\newblock {\em Annals of Statistics\/}~{\em 41}, 2097--2122.

\bibitem[\protect\citeauthoryear{Bickel, Choi, Chang, and Zhang}{Bickel, Choi,
  Chang and Zhang}{2013}]{BickelChoiChangZhang2013}
Bickel, P.~J., D.~Choi, X.~Chang, and H.~Zhang (2013).
\newblock Asymptotic normality of maximum likelihood and its variational
  approximation for stochastic blockmodels.
\newblock {\em Annals of Statistics\/}~{\em 41}, 1922--1943.

\bibitem[\protect\citeauthoryear{Bonhomme, Jochmans, and Robin}{Bonhomme,
  Jochmans and Robin}{2016}]{BonhommeJochmansRobin2014}
Bonhomme, S., K.~Jochmans, and J.-M. Robin (2016).
\newblock Nonparametric estimation of f\mbox{}inite mixtures from repeated
  measurements.
\newblock {\em Journal of the Royal Statistical Society, Series B\/}~{\em 78},
  211--229.

\bibitem[\protect\citeauthoryear{Cardoso and Souloumiac}{Cardoso and
  Souloumiac}{1993}]{CardosoSouloumiac1993}
Cardoso, J.-F. and A.~Souloumiac (1993).
\newblock Blind beamforming for non-{G}aussian signals.
\newblock {\em IEEE-Proceedings, F\/}~{\em 140}, 362--370.

\bibitem[\protect\citeauthoryear{Celisse, Daudin, and Pierre}{Celisse, Daudin
  and Pierre}{2012}]{CelisseDaudinPierre2012}
Celisse, A., J.~J. Daudin, and L.~Pierre (2012).
\newblock Consistency of maximum-likelihood and variational estimators in the
  stochastic block model.
\newblock {\em Electronic Journal of Statistics\/}~{\em 6}, 1847--1899.

\bibitem[\protect\citeauthoryear{Daudin, Picard, and Robin}{Daudin, Picard and
  Robin}{2008}]{DaudinPicardRobin2008}
Daudin, J.~J., F.~Picard, and S.~Robin (2008).
\newblock A mixture model for random graphs.
\newblock {\em Statistical Computing\/}~{\em 18}, 173--183.

\bibitem[\protect\citeauthoryear{Erd\H{o}s and R\'enyi}{Erd\H{o}s and
  R\'enyi}{1959}]{ErdosRenyi1959}
Erd\H{o}s, P. and A.~R\'enyi (1959).
\newblock On random graphs.
\newblock {\em Publicationes Mathematicae\/}~{\em 6}, 290--297.

\bibitem[\protect\citeauthoryear{Graham, Niu, and Powell}{Graham, Niu and
  Powell}{2019}]{GrahamNiuPowell2019}
Graham, B.~S., F.~Niu, and J.~L. Powell (2019).
\newblock Kernel density estimation for undirected dyadic data.
\newblock Mimeo.

\bibitem[\protect\citeauthoryear{Hoff, Raftery, and Handcock}{Hoff, Raftery and
  Handcock}{2002}]{HoffRafteryHandcock2002}
Hoff, P., A.~E. Raftery, and M.~S. Handcock (2002).
\newblock Latent space approaches to social network analysis.
\newblock {\em Journal of the American Statistical Association\/}~{\em 97},
  1090--1098.

\bibitem[\protect\citeauthoryear{Holland, Laskey, and Leinhardt}{Holland,
  Laskey and Leinhardt}{1983}]{HollandLaskeyLeinhardt1983}
Holland, P., K.~Laskey, and S.~Leinhardt (1983).
\newblock Stochastic blockmodels: {F}\mbox{}irst steps.
\newblock {\em Social Networks\/}~{\em 5}, 109--137.

\bibitem[\protect\citeauthoryear{Janson and Nowicki}{Janson and
  Nowicki}{1991}]{JansonNowicki1991}
Janson, S. and K.~Nowicki (1991).
\newblock The asymptotic distributions of generalized {U}-statistics with
  applications to random graphs.
\newblock {\em Probability Theory and Related Fields\/}~{\em 90}, 341--375.

\bibitem[\protect\citeauthoryear{Jochmans}{Jochmans}{2022}]{Jochmans2021}
Jochmans, K. (2022).
\newblock Appendix to {E}stimation and inference for stochastic block models.
\newblock Mimeo.

\bibitem[\protect\citeauthoryear{Kasahara and Shimotsu}{Kasahara and
  Shimotsu}{2014}]{KasaharaShimotsu2014}
Kasahara, H. and K.~Shimotsu (2014).
\newblock Nonparametric identif\mbox{}ication and estimation of the number of
  components in multivariate mixtures.
\newblock {\em Journal of the Royal Statistical Society, Series B\/}~{\em 76},
  97--111.

\bibitem[\protect\citeauthoryear{Kwon and Mbakop}{Kwon and
  Mbakop}{2021}]{MbakopKwon2019}
Kwon, C. and E.~Mbakop (2021).
\newblock Estimation of the number of components of non-parametric multivariate
  finite mixture models.
\newblock {\em Annals of Statistics\/}~{\em 49}, 2178--2205.

\bibitem[\protect\citeauthoryear{Le and Levina}{Le and
  Levina}{2019}]{LeLevina2019}
Le, C.~M. and E.~Levina (2019).
\newblock Estimating the number of communities by spectral methods.
\newblock Mimeo.

\bibitem[\protect\citeauthoryear{Lei}{Lei}{2016}]{Lei2016}
Lei, J. (2016).
\newblock A goodness-of-fit test for stochastic block models.
\newblock {\em Annals of Statistics\/}~{\em 44}, 401--424.

\bibitem[\protect\citeauthoryear{Levine, Hunter, and Chauveau}{Levine, Hunter
  and Chauveau}{2011}]{LevineHunterChauveau2011}
Levine, M., D.~R. Hunter, and D.~Chauveau (2011).
\newblock Maximum smoothed likelihood for multivariate mixtures.
\newblock {\em Biometrika\/}~{\em 98}, 403--416.

\bibitem[\protect\citeauthoryear{Mariadassou, Robin, and Vacher}{Mariadassou,
  Robin and Vacher}{2010}]{MariadassouRobinVacher2010}
Mariadassou, M., S.~Robin, and C.~Vacher (2010).
\newblock Uncovering latent structure in valued graphs: {A} variational
  approach.
\newblock {\em Annals of Applied Statistics\/}~{\em 4}, 715--742.

\bibitem[\protect\citeauthoryear{McLachlan and Peel}{McLachlan and
  Peel}{2000}]{McLachlanPeel2000}
McLachlan, G.~J. and D.~Peel (2000).
\newblock {\em Finite Mixture Models}.
\newblock Wiley-Blackwell.

\bibitem[\protect\citeauthoryear{Nowicki and Snijders}{Nowicki and
  Snijders}{2001}]{NowickiSnijders2001}
Nowicki, K. and T.~A.~B. Snijders (2001).
\newblock Estimation and prediction for stochastic blockstructures.
\newblock {\em Journal of the American Statistical Association\/}~{\em 96},
  1077--1087.

\bibitem[\protect\citeauthoryear{Rohe, Chatterjee, and Yu}{Rohe, Chatterjee and
  Yu}{2011}]{RoheChatterjeeYu2011}
Rohe, K., S.~Chatterjee, and B.~Yu (2011).
\newblock Spectral clustering and the high-dimensional stochastic blockmodel.
\newblock {\em Annals of Statistics\/}~{\em 39}, 1878--1915.

\bibitem[\protect\citeauthoryear{Snijders and Nowicki}{Snijders and
  Nowicki}{1997}]{SnijdersNowicki1997}
Snijders, T.~A.~B. and K.~Nowicki (1997).
\newblock Estimation and prediction for stochastic blockmodels for graphs with
  latent block structure.
\newblock {\em Journal of Classification\/}~{\em 14}, 75--100.

\bibitem[\protect\citeauthoryear{Sussman, Tang, Fishkind, and Priebe}{Sussman,
  Tang, Fishkind and Priebe}{2012}]{SussmanTangFishkindPriebe2012}
Sussman, D.~L., M.~Tang, D.~E. Fishkind, and C.~E. Priebe (2012).
\newblock A consistent adjacency spectral embedding for stochastic blockmodel
  graphs.
\newblock {\em Journal of the American Statistical Association\/}~{\em 107},
  1119--1128.

\bibitem[\protect\citeauthoryear{Tang, Cape, and Priebe}{Tang, Cape and
  Priebe}{2022}]{TangCapePriebe2022}
Tang, M., J.~Cape, and C.~E. Priebe (2022).
\newblock Asymptotically efficient estimators for stochastic blockmodels: {T}he
  naive {MLE}, the rank-constrained {MLE}, and the spectral estimator.
\newblock {\em Bernoulli\/}~{\em 28}, 1049--1073.

\bibitem[\protect\citeauthoryear{Titterington}{Titterington}{1983}]{Titterington1983}
Titterington, D.~M. (1983).
\newblock Minimum distance non-parametric estimation of mixture proportions.
\newblock {\em Journal of the Royal Statistical Society, Series B\/}~{\em 45},
  37--46.

\bibitem[\protect\citeauthoryear{Wang and Bickel}{Wang and
  Bickel}{2017}]{Wangbickel2017}
Wang, Y.~X.~R. and P.~J. Bickel (2017).
\newblock Likelihood-based model selection for stochastic block models.
\newblock {\em Annals of Statistics\/}~{\em 45}, 500--528.

\bibitem[\protect\citeauthoryear{Yan, Sarkar, and Cheng}{Yan, Sarkar and
  Cheng}{2018}]{YanSarkarCheng2018}
Yan, B., P.~Sarkar, and X.~Cheng (2018).
\newblock Provable estimation of the number of blocks in block models.
\newblock {\em Proceedings of the 21st International Conference on Artificial
  Intelligence and Statistics\/}.

\end{thebibliography}

\end{document}